\documentclass{article}
\usepackage{spconf,amsmath,graphicx}
\usepackage{amssymb}
\usepackage{amsmath}
\usepackage{graphicx}
\usepackage{graphicx}
\usepackage{makecell}
\usepackage{float}
\usepackage{spconf,amsmath,graphicx}
\newcommand{\tensor}[1]{\boldsymbol{\mathcal{#1}}}

\newcommand{\mat}[1]{\mathbf{#1}}
\newcommand{\vect}[1]{\mathbf{#1}}

\title{Randomized Tensor Ring Decomposition and Its Application to Large-scale Data Reconstruction}
\name{Longhao Yuan$^{1, 2}$, Chao Li$^{2}$, Jianting Cao$^{4, 1, *}$ and Qibin Zhao$^{2, 3, *\thanks{*Corresponding authors: qibin.zhao@riken.jp, cao@sit.ac.jp}}$ }
\address{$^1$Graduate School of Engineering, Saitama Institute of Technology, Japan \and $^2$Tensor Learning Unit, RIKEN Center for Advanced Intelligence Project (AIP), Japan
\and $^3$School of Automation, Guangdong University of Technology, China \and $^4$School of Computer Science and Technology, Hangzhou Dianzi University, China\\}

\begin{document}

\maketitle
\begin{abstract}
Dimensionality reduction is an essential technique for multi-way large-scale data, i.e., tensor. Tensor ring (TR) decomposition has become popular due to its high representation ability and flexibility. However, the traditional TR decomposition algorithms suffer from high computational cost when facing large-scale data. In this paper, taking advantages of the recently proposed tensor random projection method, we propose two TR decomposition algorithms. By employing random projection on every mode of the large-scale tensor, the TR decomposition can be processed at a much smaller scale. The simulation experiment shows that the proposed algorithms are $4-25$ times faster than traditional algorithms without loss of accuracy, and our algorithms show superior performance in deep learning dataset compression and hyperspectral image reconstruction experiments compared to other randomized algorithms.
\end{abstract}


\section{Introduction}
With the development of data acquisition and storage technology, large-scale data (i.e., big data) becomes ubiquitous in many fields such as computational neuroscience, signal processing, machine learning and pattern recognition \cite{cichocki2014era}. Among these fields, large amounts of multi-dimensional data (i.e., tensors) of high dimensionality is generated. Big data is of large volume and complex, which is hard to process by traditional methods like singular value decomposition (SVD) and principal component analysis (PCA) due to their high computational complexity. Moreover, in order to fit in these algorithms, traditional methods need to do unfolding (matricization) operations to transform tensor data to matrices and vectors, which leads to adjacent structure information loss and redundant space cost \cite{shashua2005non}.

Tensor can retain the high-dimension structure of the data and prevent information loss. Tensor decomposition aims to approximate the tensor by the latent factors, thus transforming large-scale tensor data into a latent space of low-dimensionality and reduce the data dimensionality. CANDECOMP/PARAFAC (CP) decomposition \cite{faber2003recent} and Tucker decomposition \cite{tucker1966some} are the most classical and well-studied tensor decomposition models, after which tensor train (TT) decomposition \cite{oseledets2011tensor} and tensor ring (TR) decomposition \cite{zhao2016tensor} become popular because of their high compression performance in high-dimensional and large-scale tensor. TT and TR provide a natural solution for the `curse of dimensionality'. For instance, for an order-$N$ tensor, the space complexity of Tucker grows exponentially with N, while the cases of TT, TR and CP are linear with N. Although CP is a highly compact decomposition model of which the space complexity is also linear in N, it has difficulties in finding the optimal latent tensor factors \cite{zhou2012canonical}. 

Though tensor decomposition has the merit of data structure conservation and high data representation ability, when dealing with large-scale data, traditional deterministic algorithms like alternative least squares (ALS) and gradient descent (GD) are of low-efficiency due to their high computational cost and low convergence rate. Therefore, fast and efficient algorithms are of high demand to large-scale tensor decomposition.
Randomized technology is a powerful computation acceleration technique, it has been proposed and studied for decades \cite{halko2011finding,martinsson2011randomized}. Recently, randomness-based tensor decomposition has drawn people's attention. Literature \cite{zhou2014decomposition} proposes a randomized algorithm for large-scale tensors based on Tucker decomposition, it can process arbitrarily large-scale tensors with low multi-linear rank and shows robustness to various data set. A randomized least squares algorithm for CP decomposition is proposed in \cite{battaglino2018practical}, it is much faster than traditional CP least squares algorithm and can keep the high performance at the same time. Work in \cite{erichson2017randomized} provides a different randomized CP decomposition algorithm, they first find the CP decomposition of the small tensor which is generated by tensor random projection of the large-scale tensor, then the CP decomposition of the large-scale tensor is obtained by back projection of the CP decomposition of the small tensor.


Many of these randomized tensor decomposition algorithms are efficient and perform well in simulation experiments. However, to the best of our knowledge, randomized techniques have not been applied to TR decomposition, and few studies are conducted to explore the performance of randomized tensor decomposition algorithms in real world-data. Facing the fact that TR decomposition lacks fast and efficient algorithms for large-scale tensor, in this paper, we explore the effectiveness of tensor random projection method on TR decomposition. The main contribution of this paper is listed below:
\begin{itemize}
\item Based on tensor random projection method and traditional TR decomposition algorithms, we propose two randomized TR decomposition (rTRD) algorithms, which are suitable for fast and reliable tensor decomposition of large-scale data.
\item The proposed algorithms are compared with the traditional TR decomposition algorithms in the simulation experiment. Our algorithms obtain a significant advantage in computational speed against traditional algorithms without loss of accuracy.
\item The experiments on deep learning datasets and hyperspectral image (HSI) data are conducted. The proposed algorithms outperform the compared randomized tensor decomposition algorithms in data compression and reconstruction.
\end{itemize}







\section{Notations and Preliminaries}

\subsection{Notations}
The notations in \cite{kolda2009tensor} are adopted in this paper. Tensors of order-$N\geq 3$ are denoted by calligraphic letters, e.g., $\tensor{X}\in\mathbb{R}^{I_1\times I_2\times\cdots \times I_N}$. Scalars are denoted by normal lowercase letters or uppercase letters, e.g., $x, X \in\mathbb{R}$. Vectors are denoted by boldface lowercase letters, e.g., $\vect{x}\in\mathbb{R}^{I}$. Matrices are denoted by boldface capital letters, e.g., $\mat{X}\in\mathbb{R}^{I\times J}$. For simplicity, we define tensor sequence $\{\tensor{X}_{1},\tensor{X}_{2},\ldots,\tensor{X}_{N}\}$  as $\{\tensor{X}_n\}_{n=1}^N$ or $[\tensor{X}_n]$ The scalar sequence, matrix sequence and vector sequence are denoted by the same way. Moreover, we employ two types of tensor unfolding (matricization) operations in this paper. The first mode-$n$ unfolding \cite{kolda2009tensor} of tensor $\tensor{X}  \in\mathbb{R}^{I_1\times I_2\times\cdots \times I_N}$ is denoted by $\mat{X}_{(n)}\in\mathbb{R}^{I_n \times  {I_1 \cdots I_{n-1} I_{n+1} \cdots I_N}}$, and the second mode-$n$ unfolding of tensor $\tensor{X}$ which is often used in TR operations \cite{zhao2016tensor} is denoted by $\mat{X}_{<n>}\in\mathbb{R}^{I_n \times  {I_{n+1} \cdots I_{N} I_{1} \cdots I_{n-1}}}$. In addition, the Frobenius norm of $\tensor{X}$ is defined by $\left \| \tensor{X} \right \|_F=\sqrt{\langle \tensor{X},\tensor{X} \rangle}$, where $\langle \cdot,\cdot \rangle$ is the inner product operation.


\subsection{Tensor Ring Decomposition}

Tensor ring (TR) decomposition is a more general decomposition than tensor-train (TT) decomposition, and it represents a tensor with large dimension by circular multilinear products over a sequence of low dimension cores (TR factors). All of the TR factors are order-three tensors, and are denoted by $\tensor{G}_{n} \in\mathbb{R}^{R_{n} \times I_{n} \times R_{n+1}}$, $n=1,\ldots,N$. In the same way as TT, the TR decomposition linearly scales to the dimension of the tensor, thus it can overcome the `curse of dimensionality'. $R_1, R_2,\ldots,R_N$ denotes TR-rank which controls the model complexity of TR decomposition. The TR decomposition relaxes the rank constraint on the first and last core of TT to $R_1=R_{N+1}$, while the original constraint on TT is rather stringent, i.e., $R_1=R_{N+1}=1$. TR applies trace operation and all the TR factors are constrained to be third-order equivalently. In this case, TR can be considered as a linear combination of TT and thus it offers a more powerful and generalized representation ability than TT. The element-wise relation and global relation of TR decomposition and the original tensor is given by equations (\ref{tr_relation1}) and (\ref{tr_relation2}):
{\setlength\abovedisplayskip{5pt}
\setlength\belowdisplayskip{-10pt}
\begin{equation}
\label{tr_relation1}
\tensor{X}(i_1,i_2,\ldots,i_N)=\text{Trace}\left \{  \prod_{n=1}^N \mat{G}_n(i_n) \right \},
\end{equation}}

\begin{equation}
\label{tr_relation2}
\mat{X}_{<n>}=\mat{G}_{n,(2)}(\mat{G}_{\neq n,<2>})^T,
\end{equation}
where $\text{Trace}\{ \cdot \}$ is the matrix trace operator, $ \mat{G}_n(i_n)\in\mathbb{R}^{R_n\times R_{n+1}}$ is the $i_n$th mode-$2$ slice of $\tensor{G}_n$, which also can be denoted by $\tensor{G}_n(:,i_n,:)$ according to Matlab syntax. $\tensor{G}_{\neq n}\in\mathbb{R}^{R_{n+1}\times \prod_{i=1, i\neq n}^N I_i \times R_n}$ is a subchain tensor by merging all TR factors except the $n$th core tensor, i.e., $\tensor{G}_{n+1},\ldots,\tensor{G}_N, \tensor{G}_1, \ldots ,\tensor{G}_{n-1}$, see more details in \cite{zhao2018learning}.

\section{Approach}
\subsection{Tensor Random Projection}
Tensor random projection (TRP) has drawn people's attention in the very recent years, and several studies has been conducted based on CP and Tucker \cite{erichson2017randomized,zhou2014decomposition}. Similar to matrix projection, TRP method aims to process random projection at every mode of the tensor, then a much smaller subspace tensor is obtained which reserves most of the actions of the original tensor. The TRP is simply formulated as follows:
\begin{equation}
\begin{aligned}
\tensor{X}&\approx\tensor{X}\times_1\mat{Q}_1\mat{Q}_1^T\times_2\cdots\times_N\mat{Q}_N\mat{Q}_N^T\\
&\approx\tensor{P}\times_1\mat{Q}_1\times_2\cdots\times_N\mat{Q}_N,
\end{aligned}
\label{eq_a1}
\end{equation}
where $\times_n$ is the mode-$n$ tensor production, see details in \cite{kolda2009tensor}, $[\mat{Q}_n]$ are orthogonal matrices, and $\tensor{P}$ is the projected tensor. After projection, the projected tensor $\tensor{P}$ is employed to calculate the desired low-rank approximation of the original large-scale tensor. The implementation details of the TRP method are illustrated in the next subsection.

\subsection{Randomized Tensor Ring Decomposition}
The problem of finding TR decomposition is formulated as the following model:
\begin{equation}
    \min_{[\tensor{G}_n]} \Vert \tensor{X}-\Psi([\tensor{G}_n]) \Vert_F^2,
    \label{obj}
\end{equation}
where $\tensor{X}$ is the target tensor, $[\tensor{G}_n]$ are the TR factors to be solved, and $\Psi(\cdot)$ is the function which transform the TR factors into the approximated tensor. In \cite{zhao2018learning}, the model is solved by various methods like TRSVD, TRALS, TRSGD, etc. However, the SVD-based and ALS-based algorithms are of high computational cost, when facing large-scale data,  tremendous computing resource is needed. In addition, though TRSGD owns low complexity on every iteration and is suitable for large-scale computation, the convergence speed is rather slow and the performance cannot be guaranteed. Under this situation, we combine the TRP technique with the traditional TR decomposition algorithms, (e.g. TRALS and TRSVD), to make it possible for fast and reliable TR decomposition of large-scale tensor. The randomized tensor ring decomposition (rTRD) algorithms which is based on ALS (i.e., rTRALS) and SVD (i.e., rTRSVD) are illustrated by Algorithm 1.
\begin{table}[H]
\footnotesize
\begin{center}
\begin{tabular}{l}
\hline
\textbf{Algorithm 1} Randomized tensor ring decomposition (rTRD)\\
\hline
\;\ 1: \textbf{Input}: Large-scale tensor $\tensor{X}\in\mathbb{R}^{I_1\times I_2\times \ldots\times I_N}$,\\
\;\ \;\;\;\;\;\;\;\;\;\;\;\;\;\;\ projection size of every mode $[K_n]$,\\
\;\ \;\;\;\;\;\;\;\;\;\;\;\;\;\;\ and TR-rank $R_1, \ldots, R_N$.\\
\;\ 2: $\textbf{Output}$: TR factors $[\tensor{G}_n]$ of the large-scale tensor $\tensor{X}$.\\
\;\ 3: \textbf{For} $n=1,\ldots,N$\\
\;\ 4: \;\;Create matrix $\mat{M}\in\mathbb{R}^{\prod_{i=1,i\neq n}^N I_i\times K_n}$ following\\
   \;\ \;\;\;\;\; the Gaussian distribution.\\
\;\ 5: \;\;$\mat{Y}$=$\mat{X}_{(n)}\mat{M}$ \% random projection\\
\;\ 6: \;\;$[\mat{Q}_n,\sim]=\text{QR}(\mat{Y})$ \% economy QR decomposition\\
\;\ 7: \;\;$\tensor{P}\leftarrow{}\tensor{X}\times_n \mat{Q}_n^T$\\
\;\ 8: \textbf{End for}\\
\;\ 9:\;\; Obtain TR factors $[\tensor{Z}_n]$ of $\tensor{P}$ by TRALS or TRSVD \cite{zhao2016tensor} .\\
10: \textbf{For} $n=1,\ldots,N$\\
11:\;\; $\tensor{G}_n=\tensor{Z}_n\times_n \mat{Q}_n$.\\
12: \textbf{End for} \\
\hline
\end{tabular}
\end{center}
\end{table}
It should be noted that for randomized algorithms, several techniques can be applied to the projection step to improve the numerical stability of the projection, thus providing higher decomposition performance. For example, adopting structured projection matrices instead of Gaussian distribution \cite{woolfe2008fast} and applying power iterations method to update the projected tensor in order to achieve fast decay of the spectrum of the mode-$n$ unfolding of the projected tensor \cite{halko2011finding}. In our paper, we only adopt the most basic TRP in order to show the direct improvements compared to the traditional decomposition algorithms.

\section{Experiment Results}
In the experiment section, we firstly investigate the influence of the size of the projected tensor, and compare our randomized algorithms with their traditional counterparts (i.e, rTRALS vs TRALS, and rTRSVD vs TRSVD). Then we conduct experiments on two large-scale deep learning datasets for fast data compression. Finally, a hyperspectral image (HSI) is employed to test the performance of our algorithm on data reconstruction and denoising. For evaluation index, we mainly adopt relative square error (RSE) which is calculated by $\text{RSE}=\Vert  \tensor{X}-\tensor{Y} \Vert_F /\Vert \tensor{X}\Vert_F$, where $\tensor{X}$ is the target large-scale tensor and $\tensor{Y}$ is the tensor approximated by the corresponding decomposition factors. All the computations are conducted on a Mac PC with Intel Core i7 and 16GB DDR3 memory.

\begin{figure}[htbp]
\begin{center}
\includegraphics[width=1\linewidth]{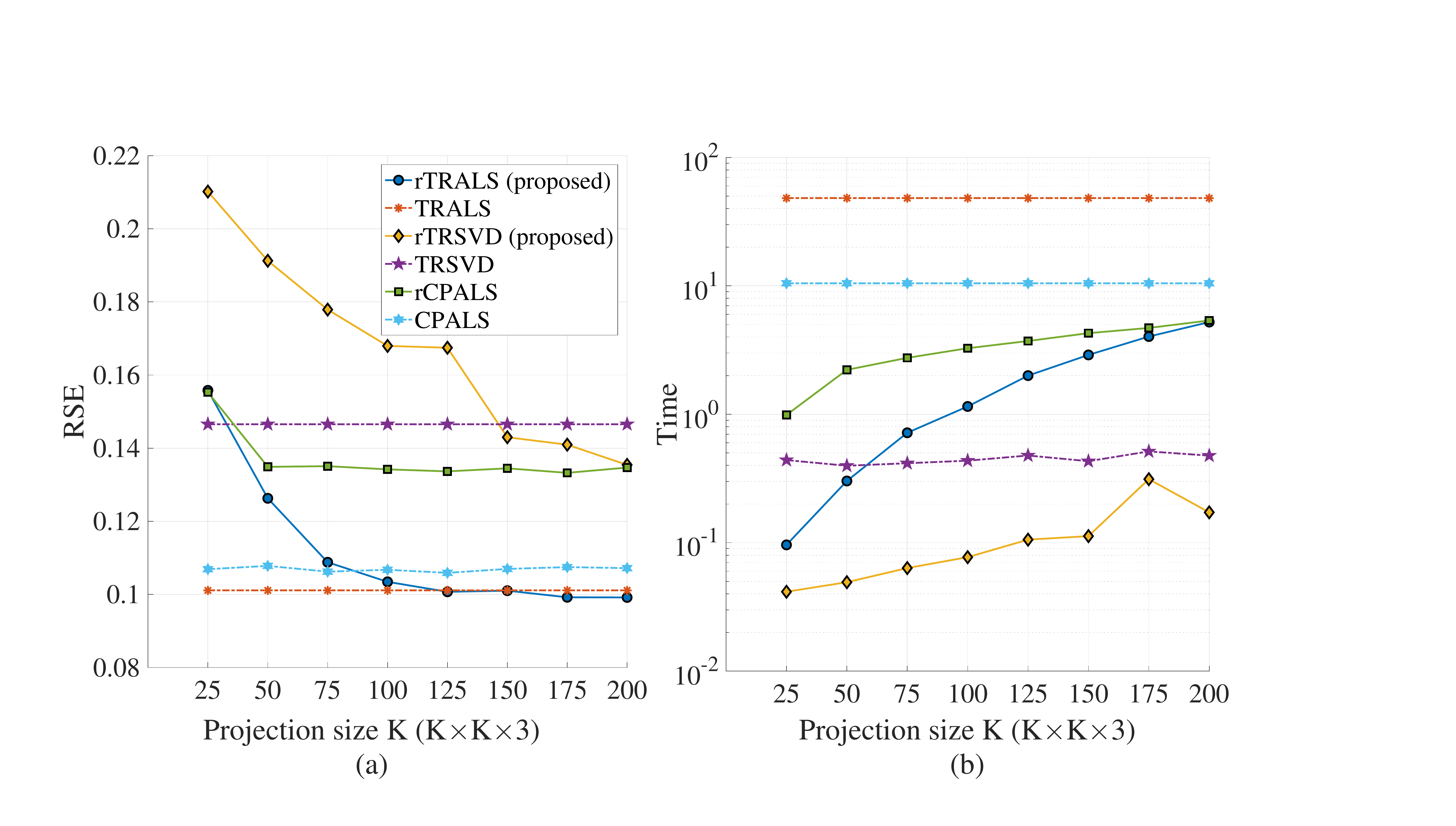}
\caption{Reconstruction results of six tensor decomposition algorithms under different tensor projection size. Figure (a) and (b) show the RSE values and the time cost respectively.}
\label{exp1}
\end{center}
\end{figure}
\setlength\tabcolsep{1.5pt} 
\begin{table*}[htb]\footnotesize
\centering
\caption{Comparison of the compression performance of randomized algorithms under two deep learning datasets. }
\label{exp2}
\setlength{\tabcolsep}{2mm}{
\begin{tabular}{c|ccc|ccc|ccc|ccc}
\hline
\hline
& \multicolumn{6}{c|}{Cifar10}& \multicolumn{6}{c}{Coil100}\\ \hline
&CR&RSE&time&CR&RSE&time&CR&RSE&time&CR&RSE&time\\
 rTRALS&102.3&0.2185 &18.29 &767.0&0.3294 &17.39 &2948.7&0.3331 &40.99 &1047.3&0.2911&42.61\\
 rTRSVD&42.64&\textbf{0.1791}&10.63&42.6&\textbf{0.1791}&10.85&175.4&\textbf{0.2669}&1.49&175.4&\textbf{0.2663}& 1.96\\
 TRSGD&102.3&0.4382&1.21e3&767.0&1.00&6.27e2&2948.7&0.4158&482.64&1047.3&0.3536&411.12\\
 rCPALS&99.0&0.2254&11.32&613.6&0.3284&10.86&3084.9&0.3434&2.12&1028.3&0.3001&5.80\\
 rTucker&100.8&0.2146&10.65&509.2&0.3058&4.61&3093.5&0.4241&0.38&1077.4&0.4680&1.98\\
\hline\hline
\end{tabular}}
\end{table*}
\subsection{Simulation}
The most important hyper-parameter of the tensor projection step is the projection size which determines the amount of residual information to be remained and controls the balance of computational speed and accuracy. In this experiment, we aim to explore how the size of the projected tensor influences the performance of our algorithms, and compare the performance with the related tensor decomposition algorithms. Except for our proposed algorithms, the rCPALS \cite{erichson2017randomized} which is the most related method is also adopted in this experiment. The counterparts of the three randomized algorithms are TRALS, TRSVD \cite{zhao2016tensor} and CPALS \cite{kolda2009tensor} respectively. We choose a RGB image of size $1024\times 1024\times 3$ as the simulation data. The projection size of order-$1$ and order-$2$ of the tensor data are chosen from $\{25,50,75,100,125,150,175,200\}$, and the order-$3$ of the tensor remains as $3$. As for parameter settings, we set the TR-rank as $\{10,10,10\}$, CP-rank as $50$, and the maximum iteration as $50$ for ALS-based algorithms. For TRSVD and rTRSVD, only one iteration is needed and the TR-rank is automatically chosen, so we only set the tolerance as 0.15. Figure \ref{exp1} shows the approximation error (RSE) and computation time of the compared algorithms. When the projection size reaches a specific value, the performance of the randomized algorithms remain steady and similar performance with their counterparts are obtained. At the steady points where the performance of the algorithm pairs are similar, from time graph we can see, rTRALS is about 24 times faster than TRALS (2.0s vs 48.1s), and rTRSVD is about 4 times faster than TRSVD (0.11s vs 0.43s).

\subsection{Deep Learning Dataset Compression}
In this section, we aim to compare the compression performance and running time of our proposed algorithms and other randomized tensor decomposition method on two deep learning datasets (i.e., CIFAR10 \cite{krizhevsky2014cifar} of size $32\times32\times 3 \times 50000$ (training data) with $1.5\times10^8$ entries, COIL100 \cite{nayar1996columbia} of size $32\times32\times 3 \times 72\times 100$ with $2.2\times10^7$ entries ). The traditional algorithms will be inefficient because the datasets are too large, so we only compare with algorithms suitable for large-scale data, i.e., TRSGD \cite{zhao2018learning}, rTucker \cite{zhou2014decomposition} and rCP \cite{erichson2017randomized}. The compression ratio CR is calculated by CR=Num/Np, where Num is the total entries of the data and Np is the number of model parameters. CR is controlled by different rank selection, and for rTRSVD, we set the tolerance as $0.15$ for automatical rank selection. Table \ref{exp2} shows the compression error and time cost of all the compared algorithms. rTRSVD and rTRALS show high accuracy and speed in all the situations, while TRSGD is much slower and obtains relatively low accuracy. Though rCPALS and rTucker are fast, the accuracy is behind our algorithms.

\subsection{Hyperspectral Image Denoising}
\begin{figure}[H]
\setlength{\abovedisplayskip}{3pt}
\setlength{\belowdisplayskip}{-100pt}
\begin{center}
\includegraphics[width=1\linewidth]{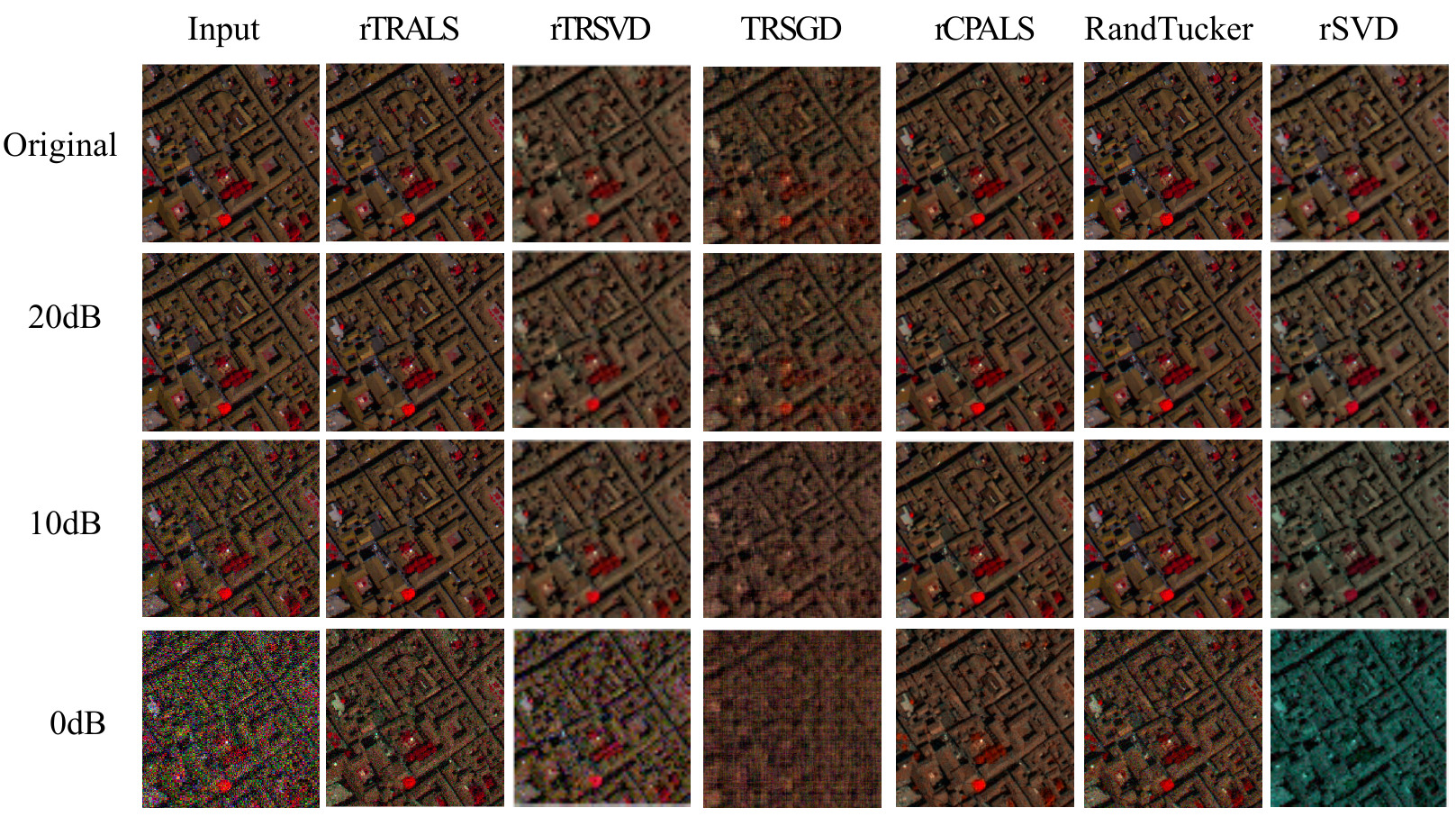}
\caption{Visual results of HSI data reconstruction with different noise}
\label{exp3}
\end{center}
\end{figure}

Hyperspectral image (HSI) is a typical type of natural order-three tensor (i.e., $height\times weights \times bands$) with large-scale. For HSI image, the spectrum-mode (mode-$3$) is usually considered to have strong low-rankness, so the projection of mode-3 can largely reduce computational cost. In this experiment, we also employ rSVD \cite{halko2011finding} which is often used in HSI image processing and rSVD is done by mode-$3$ unfolding operation. The projection size of all the algorithms are set as $100\times 100 \times 6$ for the tested $200\times 200 \times 80$ HSI image, and other parameters are set to get the best performance. Figure \ref{exp3} and Table \ref{exp3_} show the visual and numerical results respectively. rTRALS outperforms the compared algorithms in the experiment.

\vspace{-20pt}

\begin{table}[H]\footnotesize
\setlength{\abovedisplayskip}{3pt}
\setlength{\belowdisplayskip}{3pt}
\centering
\caption{Numerical results of HSI data reconstruction with different noise}
\vspace{0.1cm}
\setlength{\tabcolsep}{0.9mm}{
\begin{tabular}{c|c|c|c|c|c|c|c}
\hline
\hline
Noise&&rTR-ALS&rTR-SVD&TR-SGD&rCP-ALS&rTucker&rSVD\\
\hline
  -& \makecell[cc]{RSE \\Time} &  \makecell[cc]{\textbf{0.0150}\\ 60.01}  &   \makecell[cc]{0.149\\\textbf{0.45}}   &   \makecell[cc]{ 0.249\\9.45}  &   \makecell[cc]{0.100\\5.38}  &   \makecell[cc]{0.0110\\0.50}    & \makecell[cc]{0.0303\\1.84}  \\
  \hline
    20dB& \makecell[cc]{RSE \\Time} &  \makecell[cc]{\textbf{ 0.0294}\\60.21}   &  \makecell[cc]{ 0.143\\1.20}  &    \makecell[cc]{0.253\\206.82  }   &  \makecell[cc]{0.101\\3.97} &   \makecell[cc]{ 0.0388\\ \textbf{0.54}}    & \makecell[cc]{ 0.0594\\ 2.33} \\
    \hline
    10dB& \makecell[cc]{RSE \\Time} &   \makecell[cc]{\textbf{0.0811}\\59.61}   &    \makecell[cc]{0.113\\1.27}   &   \makecell[cc]{0.293\\210.89}   &   \makecell[cc]{0.107\\3.91}   &    \makecell[cc]{0.114\\\textbf{0.46}}&  \makecell[cc]{0.156\\2.08}    \\
    \hline
    0dB&  \makecell[cc]{RSE \\Time} &   \makecell[cc]{\textbf{0.285}\\59.05}   &   \makecell[cc]{0.328\\0.78}  &    \makecell[cc]{  0.437\\206.62}  &    \makecell[cc]{0.166\\3.95}  &       \makecell[cc]{0.367\\\textbf{0.44}}  &\makecell[cc]{0.431\\1.87}  \\
\hline\hline
\end{tabular}}
\label{exp3_}
\end{table}

\section{Conclusion}
In this paper, by tensor random projection method, we proposed rTRALS and rTRSVD algorithms for fast and reliable tensor ring decomposition. Without losing accuracy, the two algorithms perform much faster than their counterparts and outperform the other compared randomized algorithms in deep learning dataset compression and HSI image reconstruction experiments. Randomized method is a promising aspect for large-scale data processing. For future work, we will focus on further improving the performance and applying randomized algorithms to large-scale sparse and incomplete tensors.

\clearpage
\bibliographystyle{IEEEbib}
\bibliography{paper.bib}

\end{document}